\documentclass[12pt]{article}
\usepackage{amsfonts,latexsym}

\newcommand{\PP}{{\bf P}}

\newcommand{\CC}{{\bf C}}

\newcommand{\Proof}{{\it Proof\/.\quad }}
\newcommand{\F}{{\rm I\!F}}

\newcommand{\D}{{\rm I\!D}}
\newcommand{\R}{{\rm I\!R}}

\newcommand{\N}{{\rm I\!N}}

\title{ Stochastic equations with time-dependent drift driven by Levy processes}
\author{V.~P.~KURENOK$$\\
$$Department of Natural and Applied Sciences,\\
University of Wisconsin-Green Bay,\\
2420 Nicolet Drive, Green Bay, WI 54311-7001, USA}

\begin{document}

\newtheorem{Th}{Theorem}[section]
\newtheorem{Def}[Th]{Definition}
\newtheorem{Prop}[Th]{Proposition}
\newtheorem{Le}[Th]{Lemma}
\newtheorem{Cor}[Th]{Corollary}
\newtheorem{Rems}[Th]{Remarks}
\newtheorem{Rem}[Th]{Remark}
\newtheorem{Ex}[Th]{Example}
\thispagestyle{empty}
\maketitle
\vspace*{0.5cm}

\begin{abstract}
\noindent The stochastic equation $dX_t=dS_t+a(t,X_t)dt, t\ge 0,$ is considered where $S$ is a one-dimensional Levy process with the characteristic exponent $\psi(\xi), \xi\in\R$. We prove the existence of (weak) solutions for a bounded, measurable coefficient $a$ and any initial value $X_0=x_0\in\R$ when $({\mathcal Re}\psi(\xi))^{-1}=o(|\xi|^{-1})$ as $|\xi|\to\infty$. These conditions coincide with those found by Tanaka-Tsuchiya-Watanabe (1974) in the case of $a(t,x)=a(x)$. Our approach is based on Krylov's estimates for Levy processes with time-dependent drift. Some variants of those estimates are derived in this note.
\end{abstract}
{\small{\it AMS Mathematics subject classification. Primary}\hspace{15pt} 60H10, 60J60, 60J65, 60G44}

{\small{\it Keywords and phrases.}\hspace{15pt} One-dimensional Levy processes, 
time-dependent drift, Krylov's estimates, weak convergence}

\section{Introduction and preliminaries}

The goal of this note is to construct a (weak) solution of the following stochastic differential equation
\begin{equation}\label{equation1}
dX_t=dS_t+a(t,X_t)dt,\quad t\ge 0,\quad X_0=x_0\in\R,
\end{equation}
where $a:[0,\infty)\times\R\to\R$ is a measurable drift coefficient and $S$ is a one-dimensional Levy process with $S_0=0$ and the characteristic exponent $\psi(\xi), \xi\in\R$.

\quad We shall prove the existence of solutions of the equation (\ref{equation1}) for any measurable, bounded coefficient $a(t,x)$ and any Levy process $S$ satisfying the assumption
\begin{equation}\label{ourcondition}
\frac{1}{({\mathcal Re}\psi(\xi))}=o(|\xi|^{-1})\quad \mbox{ as }\quad |\xi|\to\infty.
\end{equation}
These conditions coincide exactly with those found by H.~Tanaka, M.~Tsuchiya, and S.~Watanabe \cite{TaTsWa} who considered the time-independent equation (\ref{equation1}), the case when $a(t,x)=a(x)$.

\quad In particular, the condition (\ref{ourcondition}) is satisfied when $S$ is a symmetric stable process of index $\alpha\in(1,2]$. That is, when $\psi(\xi)=|\xi|^{\alpha}$.

\quad In contrary to analytic techniques used in \cite{TaTsWa}, our approach to studying the equation (\ref{equation1}) is probabilistic and relies on using of so-called Krylov's estimates. More precisely, let $f:\R^2\to[0,\infty)$ be a measurable function and $X$ be a stochastic process of the form
\begin{equation}\label{One}
X_t=x_0+\int\limits_0^tb_sdS_s+\int\limits_0^ta_sds,\quad x_0\in \R,\quad t\ge 0,
\end{equation}
where $(b_s)$ and $(a_s)$ are two processes such that the corresponding stochastic and Lebesgue integrals are well-defined. One is interested in obtaining the estimates of the form
\begin{equation}\label{mainestimateKR}
{\bf E}\int_0^{\infty}e^{-\phi_s}\Phi_sf(s,X_s)ds\le N\|f\|_2,
\end{equation} 
where $\|f\|_{2}:=(\int_{\R^2}f^2(y)dy)^{1/2}$ and $\phi, \Phi$ are some predictable nonnegative processes.

\quad The $L_2$-estimates of the form (\ref{mainestimateKR}) are known as Krylov's estimates because N.~V.~Krylov was first who proved them for processes $X$ of diffusion type \cite{Krylov1}, the case when $S$ is a Brownian motion. They are important in the theory of stochastic differential equations as well as in their applications such as control theory, nonlinear filtering etc. Some generalizations of Krylov's estimates for diffusion processes with jumps were obtained by S.~Anulova and H.~Pragarauskas \cite{AnPr}. A.~Melnikov  \cite{Mel} derived the estimates of the form (\ref{mainestimateKR}) for some classes of semimartingales $X$. H.~Pragarauskas studied the $L_2$-estimates for Levy processes $S$ being  symmetric stable processes of index $\alpha\in(1,2)$ when $a_s=0$ \cite{Pr1}. We refer also to \cite{KU} where the case of symmetric stable processes with index $\alpha\in(1,2)$ and $a_s\neq 0$ is discussed.

\quad We shall prove here various $L_2$-estimates for processes of the form (\ref{One}) when $b=1$.


\quad We begin with some definitions. By $\D_{[0,\infty)}(\R)$ we denote the Skorokhod space, i.e. the set of all real-valued functions $x(\cdot):[0,\infty)\to \R$ with right-continuous trajectories and with finite left limits. For simplicity, we shall write $\D$ instead of $\D_{[0,\infty)}(\R)$. We will equip $\D$ with the $\sigma$-algebra ${\mathcal D}$ generated by the Skorokhod topology. Under $\D^n$ we will understand the $n$-dimensional Skorokhod space defined as $\D^n=\D\times\dots\times\D$ with the corresponding $\sigma$-algebra ${\cal D}^n$ being the direct product of $n$ one-dimensional $\sigma$-algebras ${\cal D}$.  

\quad Let $S$ be a process with $S_0=0$ defined on a complete probability space $(\Omega, {\cal F},\PP)$ and let  $\F=({\cal F}_t)$ be a filtration on $(\Omega, {\cal F},\PP)$. We use the notation $(S,\F)$ to express that $S$ is adapted to the filtration $\F$. A process $(S, \F)$ is said to be a {\it Levy process} if trajectories of $S$ belong to $\D$ and
\[
{\bf E}\left(e^{i\xi(S_t-S_s)}|{\cal F}_s\right)=e^{-(t-s)\psi(\xi)}
\]
for all $t>s\ge 0,\xi\in{\R}$ and a continuous function $\psi:\R\to\CC$. The function $\psi$ is called the {\it characteristic exponent} of the process $S$. 

\quad It is known (cf. \cite{Be}, page 13) that
\begin{equation}\label{LKh}
\psi(\xi)=ic\xi+\frac{1}{2}Q\xi^2+\int\limits_{\R\setminus\{0\}}\Bigl(1-e^{i\xi z}+i\xi z{\bf 1}_{\{|z|<1\}}\Bigr)\nu(dz),
\end{equation}
where $c\in\R$, $Q\ge 0$, and $\nu:\R\setminus\{0\}\to[0,\infty]$ is a Borel measure such that $\int(1\land z^2)\nu(dz)<\infty$.

\quad We shall use the representation $\psi(\xi)={\mathcal Re}\psi(\xi)+i{\mathcal Im}\psi(\xi)$ where the real valued functions ${\mathcal Re}\psi(\xi)$ and ${\mathcal Im}\psi(\xi)$ denote the real and imaginary part of $\psi(\xi)$, respectively. We remark also that $\psi(-\xi)$ is the characteristic exponent of the dual process $-S$ which coincides with the complex conjugate of the characteristic exponent $\psi(\xi)$ of the given process $S$. That is, $\psi(-\xi)=\overline{\psi(\xi)}={\mathcal Re}\psi(\xi)-i{\mathcal Im}\psi(\xi)$. 

\quad The measure $\nu$ in the formula (\ref{LKh}) is called the Levy measure and describes the intensity of jumps of $S$. In particular, if $\nu=0$ and $c=0$, then the Levy process $S$ is a (standard) Brownian motion process.
If $Q=0$, then $S$ is a purely jump Levy process. Because the equation (\ref{equation1}) is well-studied in the case of Brownian motion, we  shall restrict ourselves in this note to the case of purely discontinuous Levy processes ($Q=0$).
It follows then from (\ref{LKh}) that
\[
{\mathcal Re}\psi(\xi)=\int\limits_{\R\setminus\{0\}}(1-\cos{\xi z})\nu(dz)
\] 
implying 
\begin{equation}\label{nonnegative}
{\mathcal Re}\psi(\xi)\ge 0\quad \mbox{ for all }\quad \xi\in\R.
\end{equation}

\quad A stochastic process
$(X,\F)$ with trajectories in 
$\D$ , defined on a probability space $({\Omega},{\cal F},\PP)$ 
with a filtration $\F=({\cal F}_t)_{t\geq 0}$, 
is called a \emph {(weak) solution} of the equation (\ref{equation1}) with initial value $x_0\in\R$ if 
there exists a Levy process $(S,\F)$ with a given characteristic exponent $\psi$ such that
\[
X_t=x_0+S_t+\int_0^ta(s,X_s)ds,   \quad t\ge 0\quad {\bf P}\mbox{-a.s.}
\]

\quad Obviously, $S$ is a Markov process as a process having independent increaments. Hence it can be characterized in terms of its infinitesimal generator ${\mathcal L}$ defined as
\[
({\mathcal L}g)(x)=-icg^{\prime}(x)+\int\limits_{\R\setminus\{0\}}\Bigl(g(x+z)-g(x)-{\bf 1}_{\{|z|<1\}}g^{\prime}(x)z\Bigr)\nu(z)dz
\]
for any $g\in C^2$, where $C^2$ is the set of all bounded and twice continuously differentiable functions $g:\R\to\R$ (cf. \cite{Be}, page 24).  

\quad Notice finally that the use of Fourier transform can simplify calculations when working with infinitesimal operator ${\mathcal L}$. Let $g\in L_1(\R^2)$ and 
\[
\hat g(\xi_1,\xi_2):=\int\limits_{\R^2} e^{iz_1\xi_1+iz_2\xi_2}g(z_1,z_2)dz_1dz_2
\]
be the Fourier transform of $g$. Clearly, the function $\hat g(\xi_1,\xi_2)$ can be seen as the result of taking the Fourier transform from the function $g(z_1,z_2)$ first in one variable and then in another (in any order). The following facts will be used later (cf. Proposition 2.1 in \cite{KU}). 
\begin{Prop}\label{Fouriertransform} Let ${\mathcal L}$ be the infinitesimal generator of the Levy process $S$ with the characteristic exponent $\psi$. We have:
\begin{enumerate}
\item[(i)] Assume that $g\in C^2(\R)$ and ${\mathcal L}g\in L_1(\R)$. Then
\[
\widehat{({\mathcal L}g)}(\xi)=-\psi(-\xi)\hat{g}(\xi).
\]
\item[(ii)] Let $g$ be absolutely continuous on every compact subset of $\R$ and $g^{\prime}\in L_1(\R)$. Then
\[
\widehat {g^{\prime}}(\xi)=-i\xi \hat{g}(\xi).
\]
\end{enumerate}
\end{Prop}


\section{ Krylov's estimates}
\setcounter{equation}{0}

Here we shall first derive an $L_2$-estimate for solutions of a given class of  quasilinear partial  differential equtions. This estimate is then used to derive some Krylov's estimates for Levy processes with time-dependent drift. 

\quad Let $K>0$ be a constant and $f$ be a nonnegative, measurable function such that $f\in C_0^{\infty}(\R^2)$ where $C_0^{\infty}(\R^2)$ denotes the class of all infinitely many times differentiable real valued functions with compact support defined on $\R^2$. Suppose further that $S$ is a Levy process with characteristic exponent $\psi$ on a probability space $(\Omega,{\cal F},\PP)$ with filtration $\F$. By ${\mathcal I}$ we denote the class of all $\F$-predictable one-dimensional processes $(\delta_t)$ such that $|\delta_t|\le K$.

\quad For any $(t,x)\in\R^2$ and $\lambda>0$, define the value function $v(t,x)$ as
\[
v(t,x)=\sup_{\delta\in{\mathcal I}}{\bf E}\int\limits_0^{\infty}e^{-\lambda s}f(s+t, x+X^{\delta}_s)ds,
\]
where $X^{\delta}$ is a controlled process given by
\[
dX^{\delta}_t=dS_t+\delta_tdt.
\]
Then, for the value function $v$ and the process $X^{\delta}$,  the Bellman principle of optimality can be formulated as follows (cf. \cite{Krylov1}): for any $\F$-stopping time $\tau$ it holds
\[
v(t,x)=\sup_{\delta\in{\mathcal I}}{\bf E}\Bigl\{\int\limits_0^{\tau}e^{-\lambda s}f(s+t, x+X^{\delta}_s)ds+e^{-\lambda \tau}v(\tau+t,x+X^{\delta}_{\tau})\Bigr\}.
\]
Using standard arguments, one can derive from the principle above the corresponding Bellman equation ($\delta$ is a deterministic)
\[
\sup_{|\delta|\le K}\Bigl\{{v_t(t,x)+\mathcal L}v(t,x)-\lambda v(t,x)+\delta v_x(t,x)+f(t,x)\Bigr\}=0
\]
which holds a.e. in $\R^2$. Here $v_t$ and $v_x$ denote the partial derivatives of the function $v(t,x)$ in $t$ and $x$, respectively. It is not hard to see that the Bellman equation is equivalent to the equation
\begin{equation}\label{Bellmanequationrev}
v_t+{\mathcal L}v-\lambda v+K|v_x|+f=0.
\end{equation}

\begin{Le}\label{L-estimate} Assume (\ref{ourcondition}) is satisfied. Then, for all $(t,x)\in\R^2$, it holds
\begin{equation}\label{v-estimate}
v(t,x)\le N\|f\|_2,
\end{equation}
where the constant $N$ depends on $K$ and $\psi$ only.
\end{Le}
\Proof We use here the similar approach as in the proof of Lemma 3.1 in \cite{KU}. For any function $h:\R^2\to\R$ such that $h\in L_1(\R^2)$ and any $\varepsilon>0$ we define 
\[
h^{(\varepsilon)}(t,x)=\frac{1}{\varepsilon^2}\int_{\R^2}h(t,x)q\Bigl(\frac{t-s}{\varepsilon},\frac{x-y}{\varepsilon}\Bigr)dsdy
\]
to be the $\varepsilon$-convolution of $h$ with a smooth function $q$ such that $q\in C_0^{\infty}(\R^2)$ and $\int_{\R^2}q(t,x)dtdx=1$. 

For any $\varepsilon>0$, let 
\begin{equation}\label{Bellmanequationconv}
f^{(\varepsilon)}:=\lambda v^{(\varepsilon)}-v_t^{(\varepsilon)}-{\mathcal L}v^{(\varepsilon)}-K|v_x^{(\varepsilon)}|
\end{equation}
so that
\[
\Bigl(v_t^{(\varepsilon)}+{\mathcal L}v^{(\varepsilon)}-\lambda v^{(\varepsilon)}\Bigr)^2=\Bigl(K|v_x^{(\varepsilon)}|+f^{(\varepsilon)}\Bigr)^2.
\]
Obviously, $f^{(\varepsilon)}$ is square integrable and (\ref{Bellmanequationrev}) implies that $f^{(\varepsilon)}\to f$ as $\varepsilon\downarrow 0$ a.s. in $\R^2$.

It follows that
\[
\int\limits_{\R^2}\Bigl(v_t^{(\varepsilon)}(t,x)+{\mathcal L}v^{(\varepsilon)}(t,x)-\lambda v^{(\varepsilon)}(t,x)\Bigr)^2dtdx\le
\]
\[
2K^2\int\limits_{\R^2}\Bigl(|v_x^{(\varepsilon)}|(t,x)\Bigr)^2dtdx+2\int\limits_{\R^2}\Bigl(f^{(\varepsilon)}(t,x)\Bigr)^2dtdx.
\]
Now, applying Proposition \ref{Fouriertransform}, the Parseval identity and integration by parts to the last inequality yields
\begin{equation}\label{eqrev2}
\int\limits_{\R^2}|\hat{v}^{(\varepsilon)}(\zeta,\xi)|^2\Bigl([{\mathcal Re}\psi(\xi)+\lambda]^2+[{\zeta}-{\mathcal Im}\psi(\xi)]^2\Bigr)d\zeta d\xi\le 
\]
\[
2K^2\int\limits_{\R^2}|\hat{v}^{(\varepsilon)}(\zeta,\xi)|^2{\xi}^2d\zeta d\xi+2\int\limits_{\R^2} \Bigl(\hat{f}^{(\varepsilon)}(\zeta,\xi)\Bigr)^2d\zeta d\xi.
\end{equation}

Taking into account (\ref{ourcondition}) and (\ref{nonnegative}), we conclude that there exists a constant $\lambda_0>0$ such that
\begin{equation}\label{eqrev3}
[{\mathcal Re}\psi(\xi)+\lambda_0]^2+[{\zeta}-{\mathcal Im}\psi(\xi)]^2\ge[{\mathcal Re}\psi(\xi)+\lambda_0]^2\ge 4K^2{\xi}^2
\end{equation}
for all $\zeta,\xi\in\R$.

Combining the inequalities (\ref{eqrev2}) and (\ref{eqrev3}), we obtain for all $\lambda\ge \lambda_0$
\[
\frac{1}{2}\int\limits_{\R^2}|\hat{v}^{(\varepsilon)}(\zeta,\xi)|^2\Bigl([{\mathcal Re}\psi(\xi)+\lambda]^2+[{\zeta}-{\mathcal Im}\psi(\xi)]^2\Bigr)d\zeta d\xi\le
\]
\begin{equation}\label{Fourierest}
 2\int\limits_{\R^2} \Bigl(\hat{f}^{(\varepsilon)}(\zeta,\xi)\Bigr)^2d\zeta d\xi.
\end{equation}
Let
\[
N_1:=\int\limits_{\R^2}\Bigl([{\mathcal Re}\psi(\xi)+\lambda]^2+[{\zeta}-{\mathcal Im}\psi(\xi)]^2\Bigr)^{-1}d\zeta d\xi.
\]
The constant $N_1$ depends on $K$ and $\psi$ only and is finite. Indeed, 

\[
\int\limits_{\R^2}\Bigl([{\mathcal Re}\psi(\xi)+\lambda]^2+[{\zeta}-{\mathcal Im}\psi(\xi)]^2\Bigr)^{-1}d\zeta d\xi=
\]
\[
\int\limits_{\R}\Bigl\{\int\limits_{\R}\frac{1}{[{\mathcal Re}\psi(\xi)+\lambda]^2+[{\zeta}-{\mathcal Im}\psi(\xi)]^2}d\zeta\Bigr\}d\xi=
\]
\[
\pi\int\limits_{\R}\frac{1}{|\lambda+{\mathcal Re}\psi(\xi)|}d\xi<\infty,
\]
the last inequality being true because of the assumption (\ref{ourcondition}) and condition (\ref{nonnegative}).

Using the estimate (\ref{Fourierest}) and the inverse Fourier transform yields for all $t,x\in\R$ and $\lambda\ge \lambda_0$
\[
\Bigl(v^{(\varepsilon)}(t,x)\Bigr)^2\le 
\]
\[
\frac{N_1}{4\pi^2}\int\limits_{\R^2}|\hat{v}^{(\varepsilon)}(\zeta,\xi)|^2\Bigl([{\mathcal Re}\psi(\xi)+\lambda]^2+[{\zeta}-{\mathcal Im}\psi(\xi)]^2\Bigr)d\zeta d\xi  \le
\]
\[
\frac{N_1}{\pi^2}\int\limits_{\R^2} \Bigl(f^{(\varepsilon)}(s,z)\Bigr)^2dsdz.
\]

The result follows then by taking the limit $\varepsilon\to 0$ in the above inequality and using the Lebesgue dominated convergence theorem.$\Box$


\quad Now, let $X$ be a solution of the equation (\ref{equation1}) so that the assumption
\begin{equation}\label{Pragcondition}
|a(t,x)|\le K\quad\mbox{for all}\quad t,x\in\R
\end{equation}
is satisfied. We are interested in $L_2$ - estimates of the form
\begin{equation}\label{Krylovestimate}
{\bf E}\int_0^{\infty}e^{-\lambda u}f(t_0+u,x_0+X_u)du\le 
N\|f\|_{2},
\end{equation}
where $t_0,x_0\in\R$.

\begin{Th}\label{Mainestimate} Suppose $X$ is a solution of the equation (\ref{equation1}) with $X_0=0$ and the assumptions (\ref{ourcondition}) and (\ref{Pragcondition}) hold. Then, for any $t_0,x_0\in\R, \lambda\ge\lambda_0,$ and any measurable function $f:\R^2\to[0,\infty)$, it holds 
\begin{equation}\label{ourestimate}
{\bf E}\int_0^{\infty}e^{-\lambda u}f(t_0+u,x_0+X_u)du\le 
N\|f\|_{2},
\end{equation}
where the constant $N$ depends on $K$ and $\psi$ only.
\end{Th}
\Proof Assume first that $f\in C_0^{\infty}(\R^2)$ so that there is a solution $v$ of equation (\ref{Bellmanequationrev}) satisfying the inequality (\ref{v-estimate}). By taking the $\varepsilon$-convolution on both sides of (\ref{Bellmanequationrev}), we obtain 
\[
v_t^{(\varepsilon)}+{\mathcal L}v^{(\varepsilon)}-\lambda v^{(\varepsilon)}+K|v_x^{(\varepsilon)}|+f^{(\varepsilon)}\le 0.
\]

Then, for all $t_0,x_0\in\R$, applying the It\^{o}'s formula  to the expression
\[
v^{(\varepsilon)}(t_0+s,x_0+X_s)e^{-\lambda s},
\]
yields
\[
{\bf E} v^{(\varepsilon)}(t_0+s,x_0+X_s)e^{-\lambda s}- v^{(\varepsilon)}(t_0, x_0)=
\]
\[
{\bf E}\int_0^se^{-\lambda u}\Bigl[{\mathcal L}  v^{(\varepsilon)}-\lambda v^{(\varepsilon)}+a(u,X_u)v^{(\varepsilon)}_x\Bigr](t_0+u,x_0+X_u)du\le
\]
\[
{\bf E}\int_0^se^{-\lambda u}\Bigl[{\mathcal L}v^{(\varepsilon)}-\lambda v^{(\varepsilon)}+K|v^{(\varepsilon)}_x|\Bigr](t_0+u,x_0+X_u)du\le
\]
\[
-{\bf E}\int_0^se^{-\lambda u}f^{(\varepsilon)}(t_0+u,x_0+X_u)du.
\]
By Lemma \ref{L-estimate}
\[
{\bf E}\int_0^se^{-\lambda u}f^{(\varepsilon)}(t_0+u,x_0+X_u)du\le \sup_{t_0,x_0}v^{(\varepsilon)}(t_0,x_0)\le N\|f^{(\varepsilon)}\|_{2}.
\]
It remains to pass to the limit in the above inequality letting $\varepsilon\to 0$ and $s\to\infty$ and to use the Fatou's lemma.

The inequality (\ref{ourestimate}) can be extended in a standard way first to any function $f\in L_2(\R)$ and then to any non-negative, measurable function using the monotone class theorem arguments (see, for example, \cite{DelMa}, Theorem 20). $\Box$

\begin{Cor}\label{Mainestimatedrift} Suppose $X$ is a solution of the equation (\ref{equation1}) so that the assumptions (\ref{ourcondition}) and (\ref{Pragcondition}) are fulfilled. Then, for any $t\ge 0$ and any measurable function $f:\R^2\to[0,\infty)$, it holds 
\[
{\bf E}\int_0^{t}f(u,X_u)du\le 
N\|f\|_{2},
\]
where the constant $N$ depends on $K,t,$ and $\psi$ only.
\end{Cor}

\quad Now, for arbitrary but fixed $t>0, m\in\N$, define 
\[
\|f\|_{2,m,t}=\Bigl(\int_0^t\int_{-m}^m|f(s,x)|^2dxds\Bigr)^{1\over 2}
\]
to be the $L_2$-norm of $f$ on $[0,t]\times[-m,m]$. Applying (\ref{ourestimate}) to the function $\bar f(s,x)=f(s,x){\bf 1}_{[0,t]\times[-m,m]}(s,x)$, we obtain the following local version of Krylov's estimates
\begin{Cor}\label{Mainestimatelocal} Let $X$ be a solution of the equation (\ref{equation1}) with the conditions (\ref{ourcondition}) and  (\ref{Pragcondition}) being satisfied. Then, for any $t\ge 0, m\in\N,$ and any nonnegative measurable function $f$, it  follows that
\begin{equation}\label{ourestimatelocal}
{\bf E}\int_0^{t\land\tau_m(X)}f(u,X_u)du\le N\|f\|_{2,m,t},
\end{equation}
where $\tau_m(X):=\inf\{t\ge 0:|X_t|\ge m\}$ and $N$ is a constant depending on $K,\psi, m$, and $t$ only.
\end{Cor}

\section{Existence of solutions}
\setcounter{equation}{0}
 
\begin{Th}\label{boundedcoefficients} Assume the assumptions (\ref{Pragcondition}) and (\ref{ourcondition}) are true. Then, for any $x_0\in \R$, there exists a solution of the equation (\ref{equation1}).
\end{Th}
\Proof 
Because of $a$ being bounded, we can find a sequence  of functions $a_n(t,x), n\ge 1,$ such that they are globally Lipshitz continuous and uniformly bounded by the constant $K$. Then, for any $n=1,2,\dots$, the equation (\ref{equation1}) has a unique strong solution (see, for example, Theorem 9.1 in \cite{IkWa}). That is, for any fixed Levy process $S$ with the characteristic exponent $\psi$ defined on a probability space $(\Omega, {\cal F}, \PP)$, there exists a sequence of processes $X^n, n=1,2\dots,$ such that
\begin{equation}\label{strongsolution}
X^n_t=x_0+S_t+\int_0^ta_n(s,X^n_s)ds, \quad t\ge 0,\quad\PP\mbox{-a.s.}
\end{equation}

Let
\[
Y^n_t:=\int_0^ta_n(s,X^n_s)ds
\]
so that
\[
X^n=x_0+S+Y^n, n\ge 1.
\]

\quad Now we claim that the sequence of 3-dimensional processes
$Z^n:=(X^n,Y^n,S)$, $n\ge 1$, is tight in the sense of weak convergence in $(\D^3,{\cal D}^3)$. Due to the Aldous' criterion (\cite{Ald}), we have only to show that 
\begin{equation}\label{Aldous1}
\lim_{l\to\infty}\limsup_{n\to\infty}\PP\Bigl(\sup_{0\le s\le t}\|Z^n_s\|>l\Bigr)=0
\end{equation}
for all $t\ge 0$ and
\[
\limsup_{n\to\infty}\PP\Bigl(\|Z^n_{t\land(\tau^n+r_n)}-Z^n_{t\land\tau^n}\|>\varepsilon\Bigr)=0
\]
for all $t\ge 0, \varepsilon>0$, every sequence of $\F$-stopping times $\tau^n$, and every sequence of real numbers $r_n$ such that $r_n\downarrow 0$. We use $\|\cdot\|$ to denote the Euclidean norm of a vector.

\quad It is obvious that both conditions are satisfied because of the uniform boundness of the coefficients $a_n, n\ge 1$. 

\quad Since the sequence $\{Z^n \}$ is tight, there exists a subsequence $\{n_k\},k=1,2,\dots$, a probability space $(\bar\Omega,\bar{\cal F},\bar\PP)$ and the process $\bar Z$ on it with values in $(\D^3,{\cal D}^3)$ such that $Z^{n_k}$ converges weakly (in distribution) to the process $\bar Z$ as $k\to\infty$. For simplicity, let $\{n_k\}=\{n\}$.  

\quad According to the embedding principle of Skorokhod (see, e.g. Theorem 2.7 in \cite{IkWa}), there exists a probability space $(\tilde\Omega,\tilde{\cal F},\tilde\PP)$ and the processes $\tilde Z=(\tilde X,\tilde Y,\tilde S),\quad \tilde Z^n=(\tilde X^n,\tilde Y^n,\tilde S^n),\quad n=1,2,\dots,$ on it such that 
\begin{enumerate}
\item[i)] $\tilde Z^n\to\tilde Z$ as $n\to\infty$ $\tilde\PP$-a.s.
\item[ii)] $\tilde Z^n=Z^n$ in distribution for all $n=1,2,\dots.$
\end{enumerate}
Using standard measurability arguments (\cite{Krylov1}, chapter 2), one can prove that the processes $\tilde S^n$ and $\tilde S$ are Levy processes with the characteristic exponent $\psi$ with respect to the augmented filtrations $\tilde\F^n$ and $\tilde\F$ generated by processes $\tilde Z^n$ and $\tilde Z$, respectively. 

\quad Using the properties i), ii), and the equation (\ref{strongsolution}), one can show (cf. \cite{Krylov1}, chapter 2) that
\[
\tilde X^n_t=x_0+\tilde S^n_t+\int_0^ta_n(s,\tilde X^n_s)ds, \quad t\ge 0,\quad\tilde\PP\mbox{-a.s.}
\]

On the other hand,  the same properties and the quasi-left continuity of the the processes $\tilde X^n$ yield
\begin{equation}\label{convergence1}
\lim_{n\to\infty}\tilde X^n_t=\tilde X_t, \quad t\ge 0,\quad \tilde\PP\mbox{-a.s.}
\end{equation}
Therefore, in order to show that the process $\tilde X$ is a solution of the equation (\ref{equation1}), it suffices to verify that, for all $t\ge 0$,
\begin{equation}\label{convergence3}
\lim_{n\to\infty}\int_0^ta_n(s,\tilde X^n_s)ds=\int_0^ta(s,\tilde X_s)ds\quad \tilde\PP\mbox{- a.s. }
\end{equation}

The following fact can be proven similar as Lemma 4.2 in \cite{KU}.
\begin{Le}\label{limitprocess} For any Borel measurable function $f:\R^2\to[0,\infty)$ and any $t\ge 0,$ there exists a sequence $m_k\in (0,\infty), k=1,2,\dots$ such that $m_k\uparrow\infty$ as $k\to\infty$ and it holds
\[
{\bf \tilde E}\int_0^{t\land\tau_{m_k}(\tilde X)}f(s,\tilde X_s)ds\le N \|f\|_{2,m_k,t},
\]
where the constant $N$ depends on $K,\psi, t$ and $m_k$ only.
\end{Le}

\quad Without loss of generality, we can assume in the lemma above that $\{m_k\}=\{m\}$. Now, to prove (\ref{convergence3}), it is enough to verify that for all $t\ge 0$ and $\varepsilon>0$ we have 
\begin{equation}\label{convergence3'}
\lim_{n\to\infty}\tilde\PP\Bigl(|\int_0^{t}a_n(s,\tilde X^n_s)ds-\int_0^{t}a(s,\tilde X_s)ds|>\varepsilon\Bigr)=0.
\end{equation}
\quad In order to prove (\ref{convergence3'}) we estimate for a fixed $k\in\N$
\[
\tilde \PP\Bigl(|\int_0^{t}a_n(s,\tilde X^n_{s})ds-\int_0^{t}a(s,\tilde X_{s})ds|>\varepsilon\Bigr)\le
\]
\[
\tilde \PP\Bigl(|\int_0^{t}a_{k}(s,\tilde X^n_{s})ds-\int_0^{t}a_{k}(s,\tilde X_{s})ds|>{\varepsilon\over 3}\Bigr)
\]
\[
+ \tilde \PP\Bigl(|\int_0^{t\land\tau_m(\tilde X^n)}[a_{k}-a_n](s,\tilde X^n_{s})ds|>{\varepsilon\over 3}\Bigr)
\]
\[
+ \tilde \PP\Bigl(|\int_0^{t\land\tau_m(\tilde X)}[a_{k}-a](s,\tilde X_{s})ds|>{\varepsilon\over 3}\Bigr)+\tilde\PP\Bigl(\tau_m(\tilde X^n)<t\Bigr)+\tilde\PP\Bigl(\tau_m(\tilde X)<t\Bigr)=
\]
\[
\Delta^1_{n,k}+\Delta^2_{n,k,m}+\Delta^3_{k,m}+\tilde\PP\Bigl(\tau_m(\tilde X^n)<t\Bigr)+\tilde\PP\Bigl(\tau_m(\tilde X)<t\Bigr).
\]
By Chebyshev's inequality and Lebesgue bounded convergence theorem, $\Delta^1_{n,k}\to 0$ as $n\to\infty$. To show that $\Delta^2_{n,k,m}\to 0$ as $n\to\infty$ and $\Delta^3_{k,m}\to 0$ as $k\to\infty$, we use first the Chebyshev's inequality and then Corollary \ref{Mainestimatelocal} and Lemma \ref{limitprocess}, respectively, to estimate

\begin{equation}\label{est4}
\Delta^2_{n,k,m}\le {3\over \varepsilon}N\|a_{k}-a_n\|_{2,m,t}
\end{equation}
and
\begin{equation}\label{est5}
\Delta^3_{k,m}\le {3\over \varepsilon}N\|a_{k}-a\|_{2,m,t}
\end{equation}
where the constant $N$ depends on $K,m,t,$ and $\psi$ only. Obviously, $\|a_n-a\|_{2,m,t}\to 0$ as $n\to\infty$ implying that the right-hand sides in (\ref{est4}) and (\ref{est5}) converge to $0$ by letting first $n\to\infty$ and then $k\to\infty$. 

\quad Because of the property $\tau_m(\tilde X^n)\to\tau_m(\tilde X)$ as $n\to\infty$ $\tilde \PP$-a.s., 
\[
\tilde\PP\Bigl(\tau_m(\tilde X^n)<t\Bigr)\to \tilde\PP\Bigl(\tau_m(\tilde X)<t\Bigr)\quad \mbox{as}\quad n\to\infty
\]
for all $m\in\N, t>0$.
Therefore, the last two terms can be made arbitrarly small by choosing large enough $m$ for all $n$ due to the fact that the sequence of processes $\tilde X^n$ satisfies the property (\ref{Aldous1}).
This proves (\ref{convergence3'}). Hence $\tilde X$ is a solution of the equation (\ref{equation1}). $\Box$

\quad If $S$ is a symmetric stable process of index $\alpha$, then $\psi(\xi)=|\xi|^{\alpha}$ so that the assumption (\ref{ourcondition}) is satisfied for all $\alpha\in(1,2)$. It amounts us to state the following

\begin{Cor} Let $S$ be a symmetric stable process of index $\alpha\in(1,2)$ and $a(t,x)$ be mesurable and bounded. Then, for any initial value $x_0\in\R$, there exists a solution of the equation (\ref{equation1}).
\end{Cor}

\vspace{0.5cm}

{\bf Acknowledgement} The author would like to thank Henrikas Pragarauskas for valuable discussions on the subject of Krylov's estimates for stable processes.

\end{document}